# A SCHEME FOR SIMULATING ONE-DIMENSIONAL DIFFUSION PROCESSES WITH DISCONTINUOUS COEFFICIENTS

By Antoine Lejay[1] and Miguel Martinez

*Projet OMEGA, INRIA*

The aim of this article is to provide a scheme for simulating diffusion processes evolving in one-dimensional discontinuous media. This scheme does not rely on smoothing the coefficients that appear in the infinitesimal generator of the diffusion processes, but uses instead an exact description of the behavior of their trajectories when they reach the points of discontinuity. This description is supplied with the local comparison of the trajectories of the diffusion processes with those of a skew Brownian motion.

**1. Introduction.** The aim of this article is to provide a scheme for the simulation of the one-dimensional diffusion process $X$ generated by the differential operator

$$L = \frac{\rho}{2} \frac{d}{dx}\left(a \frac{d}{dx}\right) + b \frac{d}{dx}. \tag{1}$$

Here, $a$, $\rho$ and $b$ denote piecewise smooth functions that may have an infinite number of discontinuities on a countable set of points $\mathcal{J}$. We assume, however, that $\mathcal{J}$ has no cluster points and that the functions $a$, $\rho$ and $b$ have everywhere left and right limits. The triplet $(a, \rho, b)$ will be called the *characteristic* of $L$.

If $a$ belongs to $\mathcal{C}^1$, $L$ may be transformed into

$$L = \frac{a\rho}{2} \frac{d^2}{dx^2} + \left(\frac{a'}{2} + b\right)\frac{d}{dx},$$

and we see that, even if the coefficients are smooth, using a Euler scheme for the simulation of $X$ requires us to compute the derivative of $a$, which may

Received September 2004; revised June 2005.
[1]Supported in part by the "Groupe de Recherche MOMAS" funded by ANDRA, BRGM, CEA, CNRS and EDF.

be quite expensive from a numerical point of view. However, if $a = 1$, it has been proved in [39] that the Euler scheme converges but, yet, specifying its speed of convergence still remains an open and challenging problem.

The scheme presented in this article is different from the Euler scheme. It should rather be seen as a variation of the well-known random walk on spheres.

The basic idea is the following. First, we replace the differential operator $L$ by another one whose coefficients are piecewise constant, which provides good approximations of the solutions of the elliptic and parabolic PDEs involving $L$ (see Section 8 for a computation of the error). Second, we simulate the stochastic process generated by the approximation of $L$ at a given time using a description of its behavior when it is around a point where $a$ or $\rho$ (or both) are discontinuous. Mainly, we compute quantities related to the first exit time and position on some intervals. This method is exact in the case of piecewise constant coefficients when $b = 0$ because it describes correctly the behavior of the diffusion process when it reaches a point in $\mathcal{J}$.

Let us explain our approach with a simple example: assume that $b = 0$ and suppose that $(a(x), \rho(x)) = (a^+, \rho^+)$ if $x \geq 0$ and that $(a(x), \rho(x)) = (a^-, \rho^-)$ if $x < 0$. Here, $a^\pm$, $\rho^\pm$ are positive constants. Let $u$ be the weak solution of the parabolic problem associated with $L$:

(2) $$\frac{\partial u(t,x)}{\partial t} = Lu(t,x) \quad \text{and} \quad u(0,x) = f(x).$$

It is well known that this problem is equivalent to the following *transmission problem* (see [18], e.g.):

(3) $$\begin{cases} \dfrac{\partial u(t,x)}{\partial t} = \dfrac{a(x)\rho(x)}{2} \triangle u(t,x), & \text{on } \mathbb{R}_+^* \times \mathbb{R}_+^* \text{ and on } \mathbb{R}_+^* \times \mathbb{R}_-^*, \\ a^+ \nabla u(t,0+) = a^- \nabla u(t,0-), & \text{for } t > 0. \end{cases}$$

Now, let us introduce $\Phi$ defined as follows:

$$\Phi(x) = \frac{x}{\sqrt{a^+ \rho^+}} \mathbf{1}_{\{x \geq 0\}} + \frac{x}{\sqrt{a^- \rho^-}} \mathbf{1}_{\{x \leq 0\}}.$$

Since the function $u$ is of class $\mathcal{C}^2$ on $\mathbb{R}_+ \times \mathbb{R} \setminus \{0\}$, it is easy to check that $v(t,x) := u(t, \Phi(x))$ is the solution of another transmission problem

(4) $$\begin{cases} \dfrac{\partial v(t,x)}{\partial t} = \dfrac{1}{2} \triangle v(t,x), & \text{on } \mathbb{R}_+^* \times \mathbb{R}_+^* \text{ and on } \mathbb{R}_+^* \times \mathbb{R}_-^*, \\ \dfrac{\sqrt{a^+}}{\sqrt{\rho^+}} \nabla v(t,0+) = \dfrac{\sqrt{a^-}}{\sqrt{\rho^-}} \nabla v(t,0-), & \text{for } t > 0. \end{cases}$$

Formally, this new transmission problem (4) is also equivalent to the parabolic problem $\frac{\partial v}{\partial t} = \widehat{L} v(t,x)$, where $\widehat{L}$ is the formal differentiable operator $\frac{1}{2} \triangle +$

$\beta\delta_0\nabla$, with

$$\beta = \frac{\sqrt{a(0+)/\rho(0+)} - \sqrt{a(0-)/\rho(0)-}}{\sqrt{a(0+)/\rho(0+)} + \sqrt{a(0-)/\rho(0-)}} \in (-1,1).$$

As shown in [30, 31] (see also [11]), this differentiable operator $\widehat{L}$ is the infinitesimal generator of the skew Brownian motion. This process can be constructed from a reflected Brownian motion by simply choosing the sign of each excursion by tossing an independent Bernoulli random variable of parameter $(\beta+1)/2$. Heuristically, this means that the particle chooses to go on $\mathbb{R}_+$ (resp. $\mathbb{R}_-$) with probability $(\beta+1)/2$ [resp. $(1-\beta)/2$] when it reaches 0. Unfortunately, this description is not relevant due to the fact that there are infinitely many small excursions. Besides, this approach does not permit us to understand the real behavior of the particle as a function of $a^\pm$ and $\rho^\pm$. However, numerical simulation of the skew Brownian motion is easy, and using a simple and deterministic change of scale, it is possible to solve (2) using the probabilistic representation $u(t,x) = \mathbb{E}_x[f(\Phi^{-1}(Z_t^\beta))]$, where $Z^\beta$ is a skew Brownian motion of parameter $\beta$.

Instead of working with PDEs, one may prefer to use the Itô–Tanaka formula as in [29]. This can be achieved using a precise description of the process $X$ related to $L$. As shown in [20], Chapter 5, $X$ solves the following SDE:

$$X_t = X_0 + \int_0^t \sqrt{a(X_s)\rho(X_s)}\,dB_s + \frac{a(0+) - a(0-)}{a(0+) + a(0-)} L_t^0(X),$$

where $B$ is a Brownian motion and $L^0$ is the symmetric local time of $X$ at 0. The diffusion process $\Phi(X)$ is then a skew Brownian motion with the coefficient $\beta$ given above.

The basic idea of this paper is to use this kind of description because of the particular properties of the skew Brownian motion.

Although this is not the first use of the skew Brownian motion in Monte Carlo methods (see [21, 40] or, more recently, for financial applications, [9]) or in modeling (see [7] for application in ecology), there has been neither a systematic study of this process in this framework, nor a complete exposition of the interplay between the different coefficients.

The numerical method presented in this article follows the idea of [22], where diffusions with constant coefficients on the edges of a graph were simulated. As the infinitesimal generator provides locally the behavior of the particle, some properties of the skew Brownian motion allow, in this particular framework, to describe what happened at the nodes of the graph. These are the main reasons which explain our choice of approaching nonconstant coefficients using piecewise constant approximations rather than smoothing the discontinuities around the points of $\mathcal{J}$. Indeed, this last procedure turns out to be unstable in practice and very expensive numerically.

We will also provide results concerning the diffusion process generated by $L$ under more general assumptions than in [20], Chapter 5, mainly by describing it as the solution of some SDE involving the local time and being aware of the boundary conditions. In his Ph.D. thesis [27], one of the authors gives other results on diffusion processes generated by divergence-form operators, and uses also space and time transforms as a very natural tool. Among his results, he studies the speed of convergence of the Euler scheme with discontinuous coefficients obtained after making use of some other scale transform.

One could also think of using a random walk as proposed in [23] (after a proper change of the scale). The advantage of this method is that the time step is incremented with a constant value and not with a random variable. This perspective is studied by P. Étoré in the recent article [12]. We also argue that our algorithm can be implemented locally around the points where discontinuities hold. One may use more efficient algorithms (Euler scheme, Milstein scheme, . . . ) in the regions where the coefficients are smooth.

*Outline of the article.* In Sections 2 and 3 we show how to construct the stochastic process generated by $L$. In Sections 4 and 5 we study the properties of the semi-group generated by $L$ and the convergence of the solutions of the PDEs with respect to a family of differential operators. In Section 6 and 7 we study the process $X$ generated by $L$ and we show how to transform it into a process that behaves locally like a skew Brownian motion. The algorithm is explicitly given in Section 9 and its error is studied in Section 8. Finally, as an example, we show numerical results concerning the density of a doubly skewed Brownian motion and we compare this scheme with others for a differential operator with nonconstant coefficients.

HYPOTHESES. Let $\ell < r$ be two numbers that belong to $\overline{\mathbb{R}}$. It is possible that $\ell = -\infty$ and/or $r = +\infty$.

Assume that the coefficients $a$, $\rho$ and $b$ are defined on $[\ell, r]$ and satisfy

(5a)   $a$, $\rho$ and $b$ are measurable,

(5b)   for all $x \in [\ell, r]$, $\rho(x) \in [\lambda, \Lambda]$ and $a(x) \in [\lambda, \Lambda]$,

(5c)   for all $x \in [\ell, r], |b(x)| \leq \Lambda$

for some constants $\lambda, \Lambda > 0$.

Let $\nu(dx)$ be the measure $\rho(x)^{-1} dx$ on $(\ell, r)$ and $G$ denote the open set $(\ell, r)$. It is possible that $G = \mathbb{R}$. For technical reasons, we restrict ourselves to the case where $\ell$ and $r$ are simultaneously finite or infinite. In fact, all the results given here can be extended for $G = (\ell, +\infty)$ or $G = (\infty, r)$.

We define $\mathrm{L}^2(G)$ as the space of measurable functions on $G$ that are square integrable with respect to the Lebesgue measure. We also define $\mathrm{H}^1(G)$ [resp. $\mathrm{H}^1_0(G)$] as the completion of the space of smooth functions (resp. smooth functions with compact support) on $G$ with respect to the norm

$$\|f\|_{\mathrm{H}^1(G)} = \sqrt{\|f\|^2_{\mathrm{L}^2(G)} + \|\nabla f\|^2_{\mathrm{L}^2(G)}}.$$

Recall that, for any connected interval $G$, all functions in the Sobolev space $\mathrm{H}^1(G)$ have a continuous version. In all the following we will systematically identify a function in $\mathrm{H}^1(G)$ with its continuous version.

**2. Removing the drift.** We assume that $a$, $\rho$ and $b$ are smooth.

When $G$ is finite, choose $(L, \mathrm{Dom}(L))$ to be any of the differential operators:

$$\begin{cases} L \text{ is given by } (1), \\ \mathrm{Dom}(L) = \{f \in \mathcal{C}^2(G; \mathbb{R}) | f(r) = f(\ell) = 0\} \end{cases}$$

for Dirichlet boundary condition (b.c.) or

$$\begin{cases} L \text{ is given by } (1), \\ \mathrm{Dom}(L) = \{f \in \mathcal{C}^2(G; \mathbb{R}) | f'(r) = f'(\ell) = 0\} \end{cases}$$

for Neumann b.c.

When $G = \mathbb{R}$, let $(L, \mathrm{Dom}(L))$ be the differential operator:

$$\begin{cases} L \text{ is given by } (1), \\ \mathrm{Dom}(L) = \mathcal{C}^2_{\mathrm{b}}(\mathbb{R}; \mathbb{R}). \end{cases}$$

Because $a$, $\rho$ and $b$ are assumed to be smooth, it is well known that $(L, \mathrm{Dom}(L))$ is the infinitesimal generator of a continuous strong Markov process $(X, (\mathbb{P}_x)_{x \in G})$ (see [5], e.g.) and, in addition, this process is the solution to the stochastic differential equation (SDE)

$$dX_t = \sqrt{\rho a}(X_t) \, dB_t + \left(\frac{a'}{2} + b\right)(X_t) \, dt.$$

If $b = 0$, then $(L, \mathrm{Dom}(L))$ can be associated to a symmetric bilinear form

$$\mathcal{E}(u, v) = \tfrac{1}{2} \int a(x) u'(x) v'(x) \, dx$$

defined for $u, v$ in $\mathcal{C}^1_{\mathrm{b}}(G; \mathbb{R})$ through

(6) $$\mathcal{E}(u, v) = \int Lu(x) v(x) \rho^{-1}(x) \, dx$$

for $(u, v) \in \mathrm{Dom}(L) \times \mathcal{C}^1_{\mathrm{b}}(G; \mathbb{R})$.

Since the symmetry of $L$ will play an important role in most of the results given in this article, we would like to be able to transform (2) into an equivalent form without drift: as we will now see, this is possible thanks to the crucial fact that we are working in the one-dimensional space. Let us explain one way of removing the drift by transforming $a$ and $\rho$: let $\Psi$ be the function

$$\Psi(x) = \int_0^x h(x)\,dx \qquad \text{with } h(x) = 2\int_0^x \frac{b(y)}{\rho(y)a(y)}\,dy. \tag{7}$$

If $\Psi$ is given by (7) and is bounded, a simple calculation shows that

$$e^{-\Psi}\frac{\rho}{2}\frac{d}{dx}\left(ae^{\Psi}\frac{d}{dx}\right) = \frac{\rho}{2}\frac{d}{dx}\left(a\frac{d}{dx}\right) + b\frac{d}{dx}$$

and we see that

$$e^{-\Psi}\frac{\rho}{2}\frac{d}{dx}\left(ae^{\Psi}\frac{d}{dx}\right) \tag{8}$$

is a linear operator which is symmetric with respect to the measure $\rho(x)^{-1}e^{\Psi}\,dx$ and, thus, we can manipulate (6).

If $G$ is bounded, then $\Psi$ is bounded by a constant that depends only on $\lambda$, $\Lambda$ and the Lebesgue measure $\mathrm{Meas}(G)$ of $G$. Hence, this procedure is valid.

If $G$ is not bounded, then $\Psi$ is not bounded in general, and we cannot make this transformation without using further arguments. However, we claim that it is possible to replace $\Psi$, at least locally, with $\widehat{\Psi}$ defined, for example, by

$$\widehat{\Psi}(x) = \Psi(x) - \Psi(n) \qquad \text{if } x \in [n, n+1] \text{ for all } n \in \mathbb{Z}.$$

Thus, it is possible to repeat our argument and to transform $L$ locally in a linear operator which is symmetric with respect to the measure $\rho(x)^{-1}e^{\widehat{\Psi}(x)}\,dx$.

As a matter of fact, instead of considering $(a, \rho, b)$ as the characteristic of $L$, one should only consider $L$ related with the new characteristic $(e^{\widehat{\Psi}}a, e^{-\widehat{\Psi}}\rho, 0)$.

Note that the transform used here to remove the drift applies for nonsmooth coefficients $a$, $\rho$ and $b$ since it is always possible to use a regularization procedure and pass to the limit: see regularization results in Section 5.

**3. Existence of a stochastic process.** In dimension one, results on the existence of stochastic processes generated by $(L, \mathrm{Dom}(L))$ may be proved using (at least) three ways: using the properties of the density transition function, using the scale function and the speed measure, or using the Dirichlet form theory. We will see that all these methods lead to the same process.

3.1. *Using Dirichlet forms.* Let us assume at first that $b = 0$. This will allow us to use the theory of symmetric Dirichlet forms. We have seen in Section 2 that this hypothesis is not a restriction.

Let $\mathcal{E}$ be the bilinear form

$$\mathcal{E}(u,v) = \frac{1}{2} \int_G a \frac{du(x)}{dx} \frac{dv(x)}{dx} \, dx \qquad \text{for all } u, v \in \text{Dom}(\mathcal{E})$$

on $\text{L}^2(G; \nu(dx))$. The domain $\text{Dom}(\mathcal{E})$ is

$$\text{Dom}(\mathcal{E}) = \text{H}^1([\ell, r]; \nu(dx)) \simeq \text{H}^1([\ell, r]) \qquad \text{for the Neumann b.c.,}$$
$$\text{Dom}(\mathcal{E}) = \text{H}^1_0(G; \nu(dx)) \simeq \text{H}^1_0(G) \qquad \text{for the Dirichlet b.c.}$$

It is well known that $(\mathcal{E}, \text{Dom}(\mathcal{E}))$ is a regular, local Dirichlet form [13]. Hence, it generates a continuous strong Markov process $(X, (\mathcal{F}_t)_{t \geq 0}, (\mathbb{P}_x)_{x \in G})$. Besides, the process is conservative if $G = \mathbb{R}$ or $\text{Dom}(\mathcal{E}) = \text{H}^1([\ell, r])$ which corresponds to the Neumann b.c. See Lemma 1.6.5 of [13]: recurrence follows by applying Theorem 1.6.3 of [13] and the fact that, in both cases, $1 \in \text{Dom}(\mathcal{E})$ and $\mathcal{E}(1,1) = 0$.

Besides, if $G = (\ell, r)$ and $\text{Dom}(\mathcal{E}) = \text{H}^1_0(G)$, the Dirichlet form $(\mathcal{E}, \text{Dom}(\mathcal{E}))$ still possesses the local property: thus, we can repeat all the arguments of Example 4.5.1, page 166 in [13] and $X$ is absorbing at both end-points of $G$ (see also Theorem 4.2.2, page 154 in [13]). We are allowed to write $X_t = \mathbf{1}_{\{\tau \leq t\}} \delta + \mathbf{1}_{\{\tau > t\}} Y_t$, where $\delta$ is the "point at infinity" added to $\mathbb{R}$, $Y$ is the process generated by $(\mathcal{E}, \text{H}^1_0(\mathbb{R}))$ with $a = \rho = 1$ outside $G$, and $\tau = \inf\{t \geq 0 | Y_t \notin G\}$.

REMARK 1. As the dimension of the space is one, each point $x \in [\ell, r]$ is nonpolar and has a positive capacity, and all the statements of type "for quasi-every point of $[\ell, r]$" mean in fact "for every point of $[\ell, r]$."

REMARK 2. One could assume that the coefficients $a$ and $\rho$ are locally bounded and locally uniformly elliptic. In this case, the process is not necessarily conservative, which means that it may explode in finite time. This issue is discussed in [35].

3.2. *Using the properties of the semi-group.* To the Dirichlet form, $(\mathcal{E}, \text{Dom}(\mathcal{E}))$ may be associated a linear operator $(L, \text{Dom}(L))$ on $\text{L}^2(G; \nu(dx))$ defined through the relation

$$\mathcal{E}(u,v) = -\langle Lu, v \rangle_{\text{L}^2(G)\nu(dx)} \qquad \text{for all } (u,v) \in \text{Dom}(L) \times \text{Dom}(\mathcal{E}).$$

The operator $L$ is the one given by (1) with domain

$$\text{Dom}(L) = \{f \in \text{H}^1_0(G) | Lf \in \text{L}^2(G)\} \qquad \text{for the Dirichlet b.c.,}$$
$$\text{Dom}(L) = \{f \in \text{H}^1([\ell, r]) | Lf \in \text{L}^2(G)\} \qquad \text{for the Neumann b.c.}$$

It is possible to show that the operator $(L, \mathrm{Dom}(L))$ is the infinitesimal generator of a semi-group $(P_t)_{t>0}$. We will see in Section 4 that $P_t$ has a density transition function $p(t,x,y)$ and that some general estimates hold for it. These estimates [see (12), e.g.] are sufficient to ensure the existence of a strong Feller, continuous process $(X, (\mathcal{F}_t)_{t\geq 0}, (\mathbb{P}_x)_{x\in G})$.

3.3. *Using the scale function and the speed measure.* Using the scale function and the speed measure gives another way to define the stochastic process $(X, (\mathcal{F}_t)_{t\geq 0}, (\mathbb{P}_x)_{x\in G})$.

Let us set, for $x$ in $(\ell, r)$,

$$\text{(9)} \qquad h(x) = 2\int_0^x \frac{b(y)}{\rho(y)a(y)}\, dy, \qquad m(dx) = \frac{\exp(h(x))}{\rho(x)}\, dx,$$

$$\text{(10)} \qquad S(x) = \int_0^x \frac{\exp(-h(y))}{a(y)}\, dy \quad \text{and} \quad V(x) = \int_0^x m(x)\, dx.$$

For the Dirichlet b.c. (strictly speaking, the choice of $m$ in this case is that of an absorbing condition), we set $m(\{y\}) = +\infty$ for $y = \ell$ or $y = r$. For the Neumann b.c., $m(\{y\}) = 0$ for $y = r$ or $y = \ell$. On that topic, see, for example, [5] or [34]. We will see in Corollary 1 that this process corresponds to the one already constructed using Dirichlet forms.

Unless $a$ is smooth, the process $X$ is not a semi-martingale in general. It is a Dirichlet process. Nevertheless, some stochastic calculus for $X$ is possible using the theory of Dirichlet forms and time reversal techniques: see [13, 26, 32], for example.

**4. Properties of the semi-group.** Let $G$ be the open set $G = (\ell, r)$ or $G = \mathbb{R}$. For a function $u$ in $\mathrm{H}^1(G)$, we set $\Upsilon u(x) = u(x)$ if we are interested in the Dirichlet b.c. and $\Upsilon u(x) = \frac{du(x)}{dx}$ if we are interested in the Neumann b.c. Note that $\frac{du(x)}{dx}$ is a distribution in general and might not be well defined at all points, but in this section through abuse of notation, we will use this symbol as a notation for the distributional derivative of $u$.

Let us consider the parabolic partial differential equation (PDE)

$$\text{(11)} \qquad \begin{cases} \dfrac{\partial u(t,x)}{\partial t} = \dfrac{\rho(x)}{2}\dfrac{\partial}{\partial x}\left(a(x)\dfrac{\partial u(t,x)}{\partial x}\right) \\ \qquad\qquad + b(x)\dfrac{\partial u(t,x)}{\partial x}, & \text{on } (0,\infty)\times G, \\ \Upsilon u(t,x) = 0, & \text{on } (0,\infty)\times \{\ell, r\}, \\ u(0,x) = \varphi(x), & x \in G. \end{cases}$$

Let us consider the parabolic PDE (11). Unless $a$, $\rho$ and $b$ are sufficiently smooth, the solution $u$ is a *weak solution* that can only be chosen in the space $\mathcal{C}(0,T;\mathrm{L}^2(G)) \cap \mathrm{L}^2(0,T;\mathrm{H}_0^1(G))$ in the case of Dirichlet b.c. and in $\mathcal{C}(0,T;\mathrm{L}^2(G)) \cap \mathrm{L}^2(0,T;\mathrm{H}^1([\ell,r]))$ in the case of Neumann b.c. Such a solution exists and is unique as long as $\varphi$ belongs to $\mathrm{L}^2(G)$.

It is standard that $(L, \mathrm{Dom}(L))$ is the infinitesimal generator of a semi-group $(P_t)_{t>0}$ on $\mathrm{L}^2(G,\nu)$.

PROPOSITION 1. (i) *For any $t>0$, $P_t$ has a positive density function $p(t,x,y)$ with respect to the measure $\nu$. Besides,*

$$u(t,x) = \int_G p(t,x,y)\varphi(y)\,d\nu(y)$$

*is a version of the solution of* (11) *which is continuous with respect to $(t,x)$ on $(0,\infty) \times G$.*

(ii) *The function $(t,x,y) \mapsto p(t,x,y)$ is $(\alpha/2, \alpha, \alpha)$-Hölder continuous in $(t,x,y)$ on every compact of $(0,\infty) \times G$, where $\alpha$ depends only on $\lambda$, $\Lambda$ and the size of the chosen compact.*

PROOF. The existence of a density of $P_t$ is a classical result in the theory of PDEs. A proof of (ii) may be found, for example, in [2]. □

PROPOSITION 2. (i) *If $G = \mathbb{R}$ or $G$ is bounded and the Dirichlet b.c. are used, then there exists constants $C_1$ and $C_2$ depending only on $\lambda$, $\Lambda$ and $T$ such that*

$$(12) \qquad p(t,x,y) \leq \frac{C_1}{\sqrt{2\pi t}} \exp\left(-\frac{C_2|x-y|^2}{t}\right)$$

*for any $(t,x,y) \in [0,T] \times \mathbb{R} \times \mathbb{R}$. If $b=0$, then the estimate* (12) *is uniform in $T$. This bound is called Aronson's estimate. [Note that, however, D. Aronson proved these estimates in a general case, but it was initially proved simultaneously for operators of type $\nabla(a\nabla\cdot)$ by E. De Giorgi and J. Nash.]*

(ii) *If $G$ is bounded and the Neumann b.c. are used, then, for any $T>0$ and any $\theta > 1/2$, there exists some constants $C_1$, $C_2$ depending only on $\lambda$, $\Lambda$, $\theta$, $r$, $\ell$ and $T$ such that*

$$(13) \qquad p(t,x,y) \leq \frac{C_1}{t^\theta} \exp\left(\frac{C_2|x-y|^2}{t}\right)$$

*for all $t \in (0,T]$ and all $x,y \in [\ell,r]$.*

PROOF. (i) A proof of (12) can be found, for example, in [2] and in [37].

(ii) Using the continuous injection from $\mathrm{H}^1(G)$ into the set of continuous, bounded functions on $G$, it is clear that there exists a constant $C$ (depending only on $G$) such that $\|f\|_{\mathrm{L}^2(G)}^{2+4/\kappa} \leq C \|f\|_{\mathrm{H}^1(G)} \|f\|_{\mathrm{L}^1(G)}^{4/\nu}$ for every $f \in \mathrm{H}^1(G)$ and all $\kappa > 0$. Hence, with (5b), it follows that

$$\|f\|_{\mathrm{L}^2(G)nu}^{2+4/\kappa} \leq C(\mathcal{E}(f,f) + \|f\|_{\mathrm{L}^2(G)nu}^2)\|f\|_{\mathrm{L}^1(G,\nu(dx))}^{4/\kappa}$$

for all $f \in \mathrm{Dom}(\mathcal{E})$. This is the *Nash inequality*. It is then possible to apply Theorem 3.25 in [8] (see [27] for details) which yields that there exists two constants $K_1$ and $K_2$ depending only on $\kappa$, $r$ and $\ell$ such that

(14) $$p(t,x,y) \leq \frac{K_1}{t^{\kappa/2}} \exp(-|\alpha y - \alpha x| - tK_2|\alpha|^2),$$

where $\alpha = (y_0 - x_0)/Kt_0$ for an arbitrary constant $K > 0$, a time $t_0 \in (0,1]$ and fixed points $x_0, y_0 \in [\ell, r]$. For a choice of $K$ large enough in function of $K_2$ (hence, $K$ depends only on $\kappa$ and $G$) and applying (14) with $x = x_0$, $y = y_0$ and $t = t_0$, we get (13) on the time interval $(0,1]$. The Chapman–Kolmogorov equation can then be used to get (13) on any time interval $(0,T]$ for any $T > 0$. □

**5. Convergence results.** Let $(a^n, \rho^n, b^n)_{n \in \mathbb{N}}$ be a sequence of functions satisfying (5a)–(5c). Let $(L^n, \mathrm{Dom}(L^n))$ be the differential operator constructed previously, but with $(a, \rho, b)$ replaced by $(a^n, \rho^n, b^n)$.

PROPOSITION 3. *We assume that*
$$\frac{1}{a^n} \xrightarrow[n \to \infty]{\mathrm{L}^2(G)} \frac{1}{a}, \qquad \frac{1}{\rho^n} \xrightarrow[n \to \infty]{\mathrm{L}^2(G)} \frac{1}{\rho} \quad \text{and} \quad \frac{b^n}{a^n \rho^n} \xrightarrow[n \to \infty]{\mathrm{L}^2(G)} \frac{b}{a\rho}.$$

(i) *For any $\alpha > 0$ (and $\alpha = 0$ if $G$ is bounded), the weak solutions $u^n$ of the elliptic PDE $(\alpha - L^n)u^n = f$ with the Dirichlet (resp. Neumann) b.c. converge weakly in $\mathrm{H}^1(G)$ to the weak solution of $(\alpha - L)u = f$. The continuous version of $u^n(x)$ given by $\int_G g_\alpha^n(x,y)f(y)\rho^n(y)^{-1}\,dy$ with $g_\alpha^n(x,y) = \int_0^{+\infty} e^{-\alpha t} p^n(t,x,y)\,dt$ converges uniformly on each compact of $G$ to the continuous version of $u(x)$ given by $\int_G g_\alpha(x,y)f(y)\rho^{-1}(y)\,dy$, where $g_\alpha(x,y) = \int_0^{+\infty} e^{-\alpha t} p(t,x,y)\,dt$.*

(ii) *The weak solution of the parabolic PDE $\frac{\partial u^n(t,x)}{\partial t} = L^n u^n(t,x)$ with the initial condition $\varphi \in \mathrm{L}^2(G)$ converges weakly in $\mathrm{L}^2(0,T;\mathrm{Dom}(\mathcal{E}))$ to the weak solution of the parabolic PDE $\frac{\partial u(t,x)}{\partial t} = Lu(t,x)$ with the initial condition $\varphi \in \mathrm{L}^2(G)$. Moreover, the continuous version of $u^n(t,x)$ given by $\int_G p^n(t,x,y)\varphi(y)\rho^n(y)^{-1}\,dy$ converges uniformly on each compact of $\mathbb{R}_+^* \times G$ to the continuous version of $u(t,x)$ given by $\int_G p(t,x,y)\varphi(y)\rho(y)^{-1}\,dy$.*

PROOF. (i) In fact, solving $(\alpha - L^n)u^n = f$ in $\mathrm{H}^1(G, \nu(x))$ is equivalent to solving
$$\left(\frac{\alpha}{\rho^n} - \frac{1}{2}\frac{d}{dx}\left(a^n \frac{d}{dx}\right) - \frac{b^n}{\rho^n}\right)u^n(x) = \frac{f}{\rho^n}$$

in $\mathrm{H}^1(G) = \mathrm{H}^1(G; dx)$, with respect to the scalar product of $\mathrm{L}^2(G; dx)$. According to Proposition 5 and Theorem 17 in [41], this ensures the convergence of $u^n$ to $u$ in $\mathrm{H}^1(G)$.

The estimates on $p^n(t,x,y)$ given in Proposition 2 are uniform with respect to $n$. This ensures that $p^n(t,x,y)$ converges uniformly on each compact of $\mathbb{R}_+^* \times G^2$, at least along a subsequence. Indeed, its limit is necessarily $p(t,x,y)$. For each compact subset $G'$ of $G$, it is well known that any weakly convergent sequence in $H^1(G')$ converges strongly in $L^2(G')$. Combining all these facts allows us to assert that $u^n$ converges pointwise to $u$ (see [33] for the details).

(ii) If $\rho^n = 1$, the weak convergence of $u^n$ in $L^2(0,T;H^1(G))$ comes, for example, from Theorem 29 in [42], page 101.

Otherwise, we have first to use the fact that $p^n(t,x,y)$ converges pointwise to $p(t,x,y)$ in order to deduce that $u^n(t,x) = \int_G p^n(t,x,y)\varphi(y)\rho^n(y)^{-1}\,dy$ converges pointwise to $u(t,x) = \int_G p(t,x,y)\varphi(y)\rho(y)^{-1}\,dy$.

Moreover, it is standard that $u^n$ is uniformly bounded in $L^2(0,T;H^1(G))$. Thus, $u^n$ converges weakly in $L^2(0,T;G_0^1(G))$ to $u$. $\square$

Let $X^n$ be the process generated by $(L^n, \mathrm{Dom}(L^n))$.

PROPOSITION 4. *Set $G = \mathbb{R}$ or $G = (\ell, r)$ and assume that the Neumann b.c. are used in the latter case. Let $G' = (\ell', r') \subset G$, and $\tau$ (resp. $\tau^n$) be the first exit time from $G'$ for $X$ (resp. $X^n$). Under the hypotheses of Proposition 3, for any starting point $x$,*
$$\mathbb{P}_x \circ (X^n, \tau^n)^{-1} \underset{n\to\infty}{\longrightarrow} \mathbb{P}_x \circ (X, \tau)^{-1}$$
*with respect to the topology of $\mathcal{C}([0,T];\mathbb{R}) \times \mathbb{R}$ for all $T > 0$.*

PROOF. The convergence of $X^n$ to $X$ in finite-dimensional distribution follows from the convergence of the density transition function in (3), the estimates (12) or (13), and the Markov property (see, e.g., the proof of the corresponding result in [33]).

For the tightness of $(\mathbb{P}_x \circ (X_n)^{-1})_{n \in \mathbb{N}}$, according to the Aldous criterion [1], it is sufficient to prove that, for any sequence $(\tau^n, \delta_n)_{n \in \mathbb{N}}$ such that $\tau^n$ is a $\mathcal{F}_\cdot^n$-stopping time and $\delta_n > 0$ is deterministic and converges to 0, we have
$$\mathbb{P}_x[|X^n_{\tau^n + \delta_n} - X^n_{\tau^n}| > \eta] \underset{n\to\infty}{\longrightarrow} 0$$
for all $\eta > 0$. But, with (12) or (13),
$$\mathbb{P}_y[|X^n_{\delta_n}| > \eta] \leq \int_{G \setminus (y-\eta, y+\eta)} \frac{C_1}{\delta_n^\theta} \exp\left(-\frac{C_2|z-y|^2}{\delta_n}\right) dz,$$
where $\theta = 1/2$ (for $G = \mathbb{R}$) or $\theta > 1/2$ [for $G = (\ell, r)$ and Neumann b.c.]. Thus, for $n$ large enough so that $\eta^2/\delta_n \geq 1$,
$$\sup_{y \in G} \mathbb{P}_y[|X^n_{\delta_n}| > \eta] \leq C_1 \delta_n^{1/2-\theta} \int_{|z| \geq \eta/\sqrt{\delta_n}} \exp(-C_2 z 2)\,dz$$
$$\leq 2\frac{C_1}{C_2} \delta_n^{1/2-\theta} \exp(-C_2 \eta^2/\delta_n) \underset{n\to\infty}{\longrightarrow} 0.$$

The tightness of $(\mathbb{P}_x \circ (X^n)^{-1})_{n \in \mathbb{N}}$ follows immediately by application of the strong Markov property to $\mathbb{P}_x[|X^n_{\tau^n + \delta_n} - X^n_{\tau^n}| > \eta]$.

For the convergence of $(X^n, \tau^n)$, set $\Theta(x) = \inf\{t \geq 0 | x(t) \notin G'\}$ for any continuous function $x : \mathbb{R}_+ \to \mathbb{R}$. Possibly $\Theta(x) = +\infty$ if $x(t)$ stays in $G'$. It is easy to see that $\Theta$ is lower semi-continuous: that is, $\Theta(x) \leq \liminf_{n\to\infty} \Theta(x^n)$ if for all $T > 0$, $\sup_{t \in [0,T]} |x(t) - x^n(t)| \underset{n \to \infty}{\longrightarrow} 0$.

Let $x$ be a path such that there exists $(x^n)_{n \in \mathbb{N}}$ converging uniformly to $x$ on $[0, T]$ for all $T > 0$ and $\Theta(x) < \liminf_{n\to\infty} \Theta(x^n)$. It is easily seen that this means that $x$ remains in the closure of $G'$ between the time $\Theta(x)$ and $\liminf_{n\to\infty} \Theta(x^n)$. But one knows that this almost surely never happens for a trajectory of a one-dimensional regular diffusion process (see [34], Lemma V.46.1, page 273, e.g.). Thus, the set of discontinuities of $\Theta$ is a set of null measure for the distribution $\mathbb{P}_x$. Since $\tau^n = \Theta(X^n)$ by definition of $\Theta$, it follows that $(X^n, \tau^n)$ converges in distribution to $(X, \tau)$. $\square$

**6. On diffusions with discontinuous coefficients.** In this section we assume that, if there is a Dirichlet b.c. at $\ell$ (resp. $r$), then we extend the coefficients over $(-\infty, \ell]$ (resp. $[r, +\infty)$) with $\rho = a = 1$ and $b = 0$. This is justified by the results of Section 3.1 and Proposition 4.

If $\ell > -\infty$ or $r < +\infty$ and we are working with Dirichlet b.c., we extend the coefficients on $\mathbb{R}$ with $\rho = a = 1$ and $b = 0$.

Let $\mathcal{J}$ be a (countable) set of points of $(\ell, r)$, and assume the following hypotheses:

(15a) $\quad\quad\quad a, b$ and $\rho$ are right-continuous with left-limit,

(15b) $\quad\quad\quad\quad a, b$ and $\rho$ belong to $\mathcal{C}^1([\ell, r] \setminus \mathcal{J})$,

(15c) $\quad\quad\quad\quad\quad\quad \mathcal{J}$ is at most countable,

(15d) $\quad\quad$ there exists $\varepsilon > 0$ such that $|x - y| > \varepsilon$ for any $x, y \in \mathcal{J}$.

Let us remark that (15a)–(15d) ensure that the coefficients are of finite variation on $[\ell, r]$. For a function $f$ satisfying (15a)–(15d), we will denote $f'(x)$ the density of the part of its derivative which is absolutely continuous with respect to the Lebesgue measure.

PROPOSITION 5. (i) *We assume that $G = \mathbb{R}$. For any given Brownian motion $B$ constructed on a probability space $(\Omega, (\mathcal{F}_t)_{t \geq 0}, \mathcal{F}, \mathbb{P})$ with $(\mathcal{F}_t)_{t \geq 0}$ as its natural filtration (transformed to satisfy the usual conditions), $(L, \text{Dom}(L))$ is the infinitesimal generator of the unique strong solution $X$ to the SDE*

$$(16) \quad X_t = X_0 + \int_0^t \sqrt{\rho(X_s) a(X_s)}\, dB_s + \int_{\mathbb{R}} \nu(dx)\, dL^x_t(X)$$

*with*

$$\nu(dx) = \sum_{x \in \mathcal{J}} \beta(x)\delta_x + \left(\frac{a'(x)\rho(x)}{2} + b\right)\frac{dx}{a(x)} \tag{17}$$

*and*

$$\beta(x) = \frac{a(x+) - a(x-)}{a(x+) + a(x-)}. \tag{18}$$

(ii) *The operator* $(L, \mathrm{Dom}(L))$ *with the Neumann b.c. at* $r$ *and* $\ell$ *is the infinitesimal generator of the unique solution to* (16), *where* $\nu$ *and* $\beta$ *are defined by* (17) *and* (18) *and, in addition,* $\nu(\{r\}) = -1$ *and* $\nu(\{\ell\}) = 1$.

PROOF. The existence and the uniqueness of a strong solution to (16) follow from the results of J.-F. Le Gall in [23] (see also [3]). Thus, it is sufficient to prove that the process $X$ generated by $(L, \mathrm{Dom}(L))$ is a weak solution to (16).

(ii) *Case of Neumann boundary condition.*

We assume at first that $b = 0$.

It follows from the results on Dirichlet forms that a Revuz measure corresponds to a continuous additive functional under $\mathbb{P}_{x_0}$ for any $x_0 \in [\ell, r]$. For $x$ in $[\ell, r]$, let $(\widehat{L}^x_t(X))_{t \geq 0}$ be the continuous additive functional associated to the Dirac measure $\delta_x$ at $x$. We know that, for a Borel measurable bounded function $g$, the process $(\int_0^t g(X_s)\rho(X_s)\,ds)_{t \geq 0}$ is the unique continuous additive functional corresponding to the Revuz measure $g(x)\,dx$.

Let $g$ be the continuous version of a function $g$ in $\mathrm{H}^1([\ell, r])$. If $\mathrm{Id}: x \mapsto x$, an integration by parts leads to

$$\begin{aligned}
\mathcal{E}(\mathrm{Id}, g) &= \tfrac{1}{2}\int_\ell^r a(x)g'(x)\,dx \\
&= \tfrac{1}{2}\int_\ell^r a'(x)g(x)\,dx + \tfrac{1}{2}\sum_{y \in \mathcal{J}}(a(y+) - a(y-))g(y) \\
&\quad + \tfrac{1}{2}a(r)g(r) - \tfrac{1}{2}a(\ell)g(\ell) \\
&= \int_r^\ell g(x)\mu_N(dx),
\end{aligned} \tag{19}$$

where

$$\mu_N = \tfrac{1}{2}a'(x)\,dx + \sum_{y \in \mathcal{J}} \tfrac{1}{2}(a(y+) - a(y-))\delta_y + \tfrac{1}{2}a(r)\delta_r - \tfrac{1}{2}a(\ell)\delta_\ell.$$

On the other hand,

$$2\mathcal{E}(\mathrm{Id} \cdot f, f) - \mathcal{E}(f, \mathrm{Id}^2) = \int_\ell^r a(x)f(x)\,dx = \int_\ell^r f(x)\mu_M(dx). \tag{20}$$

Under $\mathbb{P}_{x_0}$, for any $x_0 \in [\ell, r]$, the process $X$ can be written $X_t = x_0 + M_t + N_t$, where $M$ is a martingale and $N$ is a continuous additive functional locally of zero quadratic variation. The bracket $\langle M \rangle$ of $M$ is a continuous additive functional characterized by the measure $\mu_M$ in (20) (see [13], equality (3.2.14), page 110), and is

$$\langle M \rangle_t = \int_0^t a(X_s)\rho(X_s)\,ds,$$

so that there exists a Brownian motion $B$, possibly on an enlarged probability space, such that

$$M_t = \int_0^t \sqrt{a(X_s)\rho(X_s)}\,dB_s.$$

The process $N$ is characterized by the measure $\mu_N$ in (19) (see [13], Theorem 5.1.3, page 187 and Corollary 5.4.1, page 224), and is then

$$N_t = \tfrac{1}{2}\int_0^t a'(X_s)\rho(X_s)\,ds + \tfrac{1}{2}a(r)\widehat{L}_t^r(X) - \tfrac{1}{2}a(\ell)\widehat{L}_t^\ell(X)$$
$$+ \sum_{y \in \mathcal{J}} \tfrac{1}{2}(a(y+) - a(y-))\widehat{L}_t^y(X).$$

Let $(L_t^x(X))_{t \geq 0}$ [resp. $L^{x+}(X)$, $L^{x-}(X)$] be the symmetric (resp. right, left) local time of $X$ at the point $x$. As both $\widehat{L}^x(X)$ and $L^x(X)$ are continuous additive functionals that increase only on $\{t \geq 0 | X_t = x\}$, there exists a real number $\gamma(x)$ such that $\widehat{L}^x(X) = \gamma(x) L^{x+}(X)$ (see, e.g., [34], Proposition 45.10, page 409).

The occupation time formula for the local time reads

(21) $$\int_\ell^r g(x+) L_t^{x+}(X)\,dx = \int_\ell^r g(X_s)\,d\langle M \rangle_s = \int_\ell^r a(X_s)\rho(X_s)g(X_s)\,ds.$$

On the other hand,

(22) $$\int_\ell^r g(x)\widehat{L}_t^x(X)\,dx = \int_0^t g(X_s)\rho(X_s)\,ds.$$

As (21) and (22) are true for any measurable and bounded function $g$,

$$\gamma(x) = \frac{1}{a(x+)}.$$

One knows that

$$L_t^{x+}(X) - L_t^{x-}(X) = 2\int_0^t \mathbf{1}_{\{X_s = x\}}\,dX_s.$$

Thus, for any $x \in \mathcal{J}$,

(23) $$L_t^{x+}(X) - L_t^{x-}(X) = (a(x+) - a(x-))\gamma(x) L_t^{x+}(X)$$

if $x \in \mathcal{J}$. As, by definition,
$$L_t^x(X) = \tfrac{1}{2}(L_t^{x+}(X) + L_t^{x-}(X)),$$
one gets easily the value of $\beta(x)$.

At point $\ell$, $\widehat{L}_t^\ell(X) = \gamma(\ell) L_t^{\ell+}(X)$ and $\widehat{L}_t^\ell(X) = \tfrac{1}{2}\widehat{L}_t^{\ell+}(X)$. Thus, $\frac{a(\ell)}{2}\widehat{L}^\ell(X) = L_t^\ell(X)$. A similar result holds for $L_t^{r-}(X)$. Hence, $X$ is a weak solution to (16).

If $b \neq 0$, substituting $ae^\Psi$ and $\rho e^{-\Psi}$, where $\Psi$ is defined in (7), to $a$ and $\rho$ yields immediately the result.

(i) *If $G = \mathbb{R}$.*

The previous computations can be used with a localization procedure. But it is possible to avoid using the theory of Dirichlet forms. With the Itô–Tanaka formula (see, e.g., [34], Chapter IV.45, page 102) and the results from [23], $X = S^{-1}(Y)$ is a solution to (16) if $Y$ is the strong solution of
$$Y_t = S(x) + \int_0^t e^{-\Psi(S^{-1}(Y_s))} \sqrt{\rho(S^{-1}(Y_s))/a(S^{-1}(Y_s))}\, dB_s$$
and the proof of Proposition 5 is now complete. □

The next results follow from the Itô–Tanaka formula (see [34], Chapter IV.45, page 102, e.g.) applied to $S(X_t)$, where $S$ is given by (10), and $X$ is the solution to (16).

COROLLARY 1. *The speed measure $m$ and the scale function $S$ of $(X, (\mathbb{P}_x)_{x \in G})$ are given in* (9) *and* (10).

REMARK 3. If $G = \mathbb{R}$, this result could also be proved using a smooth approximation of the coefficients and the results of [14].

**7. Approximation by diffusions with piecewise constant coefficients.** In this section we assume that $b = 0$, and we have seen in Section 2 that it is possible to transform $a$ and $\rho$ in order to remove the drift. Yet, the results of this section may easily be extended to the case $b \neq 0$.

7.1. *The SDEs satisfied by the approximations.* To simplify the notation, we assume that $\ell > -\infty$, $r = \infty$ and we set, for $n \in \mathbb{N}$, $\mathcal{J}^n = \{x_i | \ell = x_0^n < x_1^n < \cdots\}$ for some points $x_0^n, x_1^n, \ldots$. If $\ell = -\infty$ and $r = +\infty$, we have to pick a reference point on $\mathbb{R}$ and to use doubly indexed sequences.

Thus, we set, for $f = a$ and $f = \rho$,
$$f^n(x) = \sum_{k \geq 0} \mathbf{1}_{[x_i^n, x_{i+1}^n)}(x) f(\widehat{x}_i^n),$$
where $\widehat{x}_i^n$ is a point in $[x_i^n, x_{i+1}^n)$.

HYPOTHESIS 1. *For any $n \in \mathbb{N}$, the points $x_0^n < x_1^n < \cdots$ are chosen such that $\mathcal{J} \subset \mathcal{J}^n$, the minimal distance between two points of $\mathcal{J}^n$ is positive, and*

$$\|a - a^n\|_\infty + \|\rho - \rho^n\|_\infty \xrightarrow[n \to \infty]{} 0.$$

Since $(a, \rho)$ satisfies (15b) and we assume that $\mathcal{J} \subset \mathcal{J}^n$, it is clear that one may construct, at least on a compact subset of $(\ell, r)$, such a sequence $(\mathcal{J}^n)_{n \in \mathbb{N}}$. For each $n \in \mathbb{N}$, the piecewise constant coefficients $a^n$ and $\rho^n$ and the set $\mathcal{J}^n$ satisfy (5a)–(5b) and (15a)–(15d), so that the results of the previous sections apply.

Using the occupation time density formula, it follows from Proposition 5 that the diffusion $X^n$ is solution to the SDE

$$(24) \qquad X_t^n = X_0^n + \int_0^t \sqrt{a^n \rho^n(X_s^n)}\, dB_s + \sum_{k \geq 0} \widetilde{\beta}_k^n L_t^{x_k^n}(X^n),$$

where $B$ is a Brownian motion,

$$\widetilde{\beta}_k^n = \frac{a^n(x_k^n+) - a^n(x_k^n-)}{a^n(x_k^n+) + a^n(x_k^n-)} \qquad \text{if } k \neq 0,$$

and $L_t^x(X^n)$ is the symmetric local time of $X^n$ at point $x$ and time $t$. If the Neumann b.c. is used at $\ell$, then $\beta_0^n = 1$. If the Dirichlet b.c. is used at $\ell$, then we consider (24) up to $\tau = \inf\{t \geq 0 | X_t^n = \ell\}$.

Let $\Phi^n$ be the piecewise linear function

$$\Phi^n(x) = \sum_{k=0}^{k^n(x)-1} \frac{x_{k+1}^n - x_k^n}{\sqrt{a^n(x_k^n)\rho^n(x_k^n)}} + \frac{x - x_{k^n(x)}^n}{\sqrt{a^n(x_{k^n(x)}^n)\rho^n(x_{k^n(x)}^n)}},$$

where $k^n(x)$ is such that $x_{k^n(x)}^n \leq x < x_{k^n(x)+1}^n$. Then the symmetric Itô–Tanaka formula (see [34], Chapter IV.45, page 102, e.g.; see also [29] for a treatment of this case) applied to $X^n$ yields that $Y^n = \Phi^n(X^n)$ is the solution to the SDE

$$(25) \qquad Y_t^n = Y_0^n + B_t + \sum_{k \geq 0} \beta_k^n L_t^{y_k^n}(Y^n),$$

where $y_k^n = \Phi(x_k^n)$ and

$$(26) \quad \beta_k^n = \frac{\sqrt{a^n(x_k^n+)/\rho^n(x_k^n+)} - \sqrt{a^n(x_k^n-)/\rho^n(x_k^n-)}}{\sqrt{a^n(x_k^n+)/\rho^n(x_k^n+)} + \sqrt{a^n(x_k^n-)/\rho^n(x_k^n-)}} \qquad \text{if } k \neq 0.$$

Of course, if the Dirichlet b.c. is used at $\ell$, then we consider only $Y$ up to $\tau = \inf\{t \geq 0 | Y_t = \Phi(\ell)\}$. If the Neumann b.c. is used at $\ell$, then we set $\beta_0^n = 1$. The b.c. at $r$ is treated the same way.

The infinitesimal generator of $Y^n$ is $L^{Y^n} = \frac{1}{2}\triangle$ on $(\ell, +\infty) \setminus \mathcal{J}^n$, whose domain $\text{Dom}(L^{Y^n})$ is the closure of the set of continuous, bounded functions $f$ of class $\mathcal{C}_b^2$ on $(\ell, \infty) \setminus \mathcal{J}^n$ such that $f(\ell) = 0$ for the Dirichlet b.c. and $f'(\ell) = 0$ for the Neumann b.c. and

$$\frac{1+\beta_k^n}{2} f'(x_k^n +) = \frac{1-\beta_k^n}{2} f'(x_k^n -)$$

for each integer $k$ and $n$.

REMARK 4. If $\rho = 1$ and $a(x) = a_+$ on $\mathbb{R}_+$ and $a(x) = a_-$ on $\mathbb{R}_-$ with $a_+, a_- > 0$, then one derives from (25) and (26) that $\mathbb{P}_0[X_t > 0] = \beta = \sqrt{a_+}/(\sqrt{a_+} + \sqrt{a_-})$ for any $t > 0$. The geophysical community has already noticed that, in a heterogeneous media with a diffusion coefficient taking two values $a_+$ and $a_-$, the parameter $\beta$ gives the ratio of concentration of a fluid in the upper half-space (see [36], e.g.).

7.2. *The skew Brownian motion.* Equation (25) shows that the process $Y^n$ behaves around $x_k^n$ like a *skew Brownian motion* of parameter $(1+\beta_k^n)/2$. Various points of view and constructions of this process may be found in [4, 16, 17, 38].... Of course, the skew Brownian motion of parameter 1 (resp. 0) is the positively (resp. negatively) reflected Brownian motion (this is used to deal with the Neumann b.c.).

Let $Z^\theta$ be a skew Brownian motion of parameter $\theta \in [0,1]$, and set

$$\tau = \inf\{t \geq 0 | |Z_t^\theta| = \rho\} \tag{27}$$

for some $\rho > 0$. We will use the following construction of the skew Brownian motion, which may be found in [17], Problem 1, page 115: a skew Brownian motion can be constructed by flipping the excursions of a reflected Brownian motion with a probability $\theta$.

To be more precise, let $R$ be a reflected Brownian motion, and $\{(\ell^n, r^n)\}_{n \in \mathbb{N}}$ be the family of its excursions intervals. These intervals are such that $R_{\ell^n} = R_{r^n} = 0$, $R_t > 0$ on $(\ell^n, r^n)$, $\overline{\bigcup_{n \in \mathbb{N}}(\ell^n, r^n)} = \mathbb{R}_+$ and $(\ell^n, r^n) \cap (\ell^k, r^k) = \varnothing$ for $n \neq k$ (see [34], e.g., for the existence of these intervals). Let $e^n$ be the excursion attached to the interval $n$, that is, $e_t^n = R_{(t-\ell^n) \wedge (r^n - \ell^n)}$ [note that these intervals $(\ell^n, r^n)$ may not be ordered, which implies that $n$ is just understood as a label with a priori no other meaning]. To the excursion $n$, we associate an independent Bernoulli random variable $\sigma_n$ of parameter $\theta$ with value in $\{-1, 1\}$. The process $Z^\theta$ constructed by

$$Z_t^\theta = \sigma_n e_{t-\ell^n}^n \qquad \text{if } t \in [\ell^n, r^n]$$

is then the skew Brownian motion of parameter $\theta$.

The core idea of the algorithm is contained in the following lemma.

LEMMA 1. *The random variables $(\tau, Z_\tau^\theta)$ are independent. Moreover, the distribution of $\tau$ does not depend on $\theta$ (in particular, $\tau$ is equal in distribution to the first exit time of $[-\rho, \rho]$ for a Brownian motion), and $\mathbb{P}_0[Z_\tau^\theta = \rho] = \theta$.*

PROOF. This is a direct consequence of the previous construction of the skew Brownian motion. □

As we will see in the sequel, we will need to simulate a realization of $Z_t^\theta$ given $\{t < \tau\}$ under $\mathbb{P}_0$ for some $\rho > 0$, where $\tau$ is defined by (27).

LEMMA 2. *Let $R$ be a reflected Brownian motion. Then, if $0 \leq y_0 \leq y_1 \leq \rho$,*

$$(28) \qquad \mathbb{P}_0[Z_t^\theta \in [y_0, y_1]; t < \tau] = \theta \mathbb{P}_0[R_t \in [y_0, y_1]; t < \tau_\rho(R)]$$

*and if $-\rho \leq y_0 \leq y_1 \leq 0$,*

$$(29) \qquad \mathbb{P}_0[Z_t^\theta \in [y_0, y_1]; t < \tau] = (1-\theta)\mathbb{P}_0[R_t \in [-y_1, -y_0]; t < \tau_\rho(R)],$$

*where $\tau_\rho(R) = \inf\{t \geq 0 | R_t = \rho\}$.*

In other words, to simulate $Z_t^\theta$, one needs to simulate the position of a reflected Brownian motion $R$ at time $t$ given $t < \tau_\rho(R)$ [or to reflect the position of a Brownian motion at time $t$ given $t < \tau$, whose density is given by either (39) or (40)], and to use an independent Bernoulli random variable for the sign of $Z_t^\theta$.

PROOF OF LEMMA 2. For a trajectory of the skew Brownian motion $Z^\theta$, set $\{(\ell^n, r^n)\}_{n \in \mathbb{N}}$ as the excursions intervals of $Z^\theta$. For each $n \in \mathbb{N}$, set $e_t^n = Z_{(t-\ell^n) \wedge (r^n - \ell^n)}^\theta$ if $t \in [\ell^n, r^n]$, which are the excursions of $Z^\theta$. We assume that $0 \leq y_0 \leq y_1 \leq \rho$. Then

$$\mathbb{P}_0[Z_t^\theta \in [y_0, y_1]; \tau_{-\rho, \rho}(R)]$$
$$= \sum_{n \in \mathbb{N}} \mathbb{P}_0 \bigg[ |e_{t-\ell^n}^n| \in [y_0, y_1]; \mathrm{sgn}(e^n) = 1;$$
$$\tau_{-\rho, \rho}(e^n) > t - \ell^n; \sup_{k \text{ s.t. } r^k < \ell^n} \sup_{s \in [0, r^k - \ell^k]} |e^k(s)| < 1 \bigg],$$

where $\mathrm{sgn}(e^n)$ is the sign of the excursion $e^n$, and $\tau_{-\rho, \rho}(e^n)$ is the first time (possibly infinite) when the excursion $e^n$ reaches $-\rho$ or $\rho$. By construction (see Section 7.2), the sign $\mathrm{sgn}(e^n)$ of the excursion $e^n$ is a Bernoulli random variable of parameter $\theta$ which is independent from any other random variables involved in this construction. Equality (28) follows easily, and (29) is proved in a similar way. □

## 8. Error estimates.

8.1. *The elliptic case.* Set the Dirichlet problem

$$\begin{cases} Lu = 0 & \text{on } [r, \ell], \\ u(r) = u_r & \text{and} \quad (\ell) = u_\ell. \end{cases} \tag{30}$$

We study the error of the solution when $(a, \rho, b)$ is replaced by $(a^n, \rho^n, b^n)$ in (30) and when

$$\sup_{x \in G} \|(a^n, \rho^n, b^n)(x) - (a, \rho, b)(x)\| \xrightarrow[n \to \infty]{} 0.$$

PROPOSITION 6. *There exists a constant $C$ depending only on $\lambda$, $\Lambda$, $r$ and $\ell$ such that*

$$\sup_{x \in G} |u(x) - u^n(x)| \leq C|u_r - u_\ell| \sup_{x \in G} \|(a^n, \rho^n, b^n)(x) - (a, \rho, b)(x)\|.$$

PROOF. We use for $L$ the operator $L = 2^{-1} e^{-\Phi} \rho \frac{d}{dx}(e^\Phi a \frac{d}{dx})$. Let $\widehat{a}^n$ and $\widehat{\rho}^n$ be some approximations of $\widehat{a} = ae^\Phi$ and $\widehat{\rho} = \rho e^{-\Phi}$, such that $\widehat{a}^n$ and $\widehat{\rho}^n$ converge uniformly to $\widehat{a}$ and $\widehat{\rho}$.

Let $u^n$ be the solution of (30) with $L$ replaced by $L^n$, and $u$ the solution of (30).

Then, using $v^n = u^n - u \in H^1_0(G)$ as a test function, one gets

$$\int_G \widehat{a}^n \left(\frac{dv^n}{dx}\right)^2 dx = \int_G (\widehat{a}^n - \widehat{a}) \frac{du}{dx} \frac{dv^n}{dx} dx.$$

As $xy \leq \lambda x^2/2 + 2\lambda^{-1} y^2$ for all $x, y \geq 0$ and all $\lambda > 0$, it follows classically that

$$\frac{\lambda}{2} \left\| \frac{dv^n}{dx} \right\|^2_{L^2(G)} \leq 2\lambda^{-1} \|\widehat{a}^n - \widehat{a}\|_\infty \left\| \frac{du}{dx} \right\|^2_{L^2(G)}.$$

Let $\varphi$ be the linear function $\varphi(x) = (u_r - u_\ell)(x - \ell)/(r - \ell) + u_\ell$. Then, $u - \varphi$ belongs to $H^1_0(G)$, and then

$$\int_G a \left(\frac{du}{dx}\right)^2 dx = \frac{u_r - u_\ell}{r - \ell} \int_G a \frac{du}{dx} dx.$$

Hence, there exists some constant $C$, depending only on $\Lambda$, $r$, $\ell$, $u_r$ and $u_\ell$ such that

$$\int_G \left(\frac{du}{dx}\right)^2 dx \leq C \left(\frac{u_r - u_\ell}{r - \ell}\right)^2.$$

As $u^n$ and $u$ satisfy the same boundary condition, $v^n$ belongs to $\mathrm{H}_0^1(G)$, and from the Poincaré inequality,

$$\|v^n\|_{\mathrm{L}^2(G)}^2 \leq C \left\|\frac{dv^n}{dx}\right\|_{\mathrm{L}^2(G)}^2.$$

Moreover, $\mathrm{H}^1(G)$ can be continuously injected in the space of continuous functions on $G$. Thus, for some constant $C'$, $\sup_{x \in G} |v^n(x)|^2 \leq C' \|v^n\|_{\mathrm{H}^1(G)}^2$. It follows that

$$\sup_{x \in G} |u^n(x) - u(x)| \leq C'' \|\widehat{a}^n - \widehat{a}\|_\infty \left|\frac{u_r - u_\ell}{r - \ell}\right|$$

for some constant $C''$ depending only on $r$, $\ell$, $u_r$, $u_\ell$, $\lambda$ and $\Lambda$.

To conclude the proof, it remains to see that

$$|\widehat{a}^n(x) - \widehat{a}(x)| \leq C_1 \sup_{x \in G} |a^n(x) - a(x)|$$
$$+ C_2 \sup_{x \in G} |\rho^n(x) - \rho(x)| + C_3 \sup_{x \in G} |b^n(x) - b(x)|$$

for some constants $C_1$, $C_2$ and $C_3$ depending only $\lambda$, $\Lambda$, $r$ and $\ell$. $\square$

8.2. *The parabolic case.* The parabolic case is harder to deal with, and we are not able to give a full treatment of it.

Let $u$ be the solution of $\frac{\partial u}{\partial t} = Lu$ on $\mathbb{R}_+ \times G$, with the Dirichlet or Neumann b.c. on $\mathbb{R}_+ \times \{\ell, r\}$ and initial condition $\varphi(x)$. Let also $u^n$ be the solution of the similar problem where $L$ is replaced by $L^n$.

PROPOSITION 7. (i) *Assume that $\rho \equiv 1$. Then, for any finite open interval $(\ell', r') \subset G$, there exists a constant $C$ depending on $\lambda$, $\Lambda$, $\|\varphi\|_{\mathrm{L}^2}$, $\ell'$, $r'$ and $T$ such that*

$$\sup_{(t,x) \in [0,T] \times (\ell',r')} |u^n(t,x) - u(t,x)| \leq C(\|a^n - a\|_\infty + \|b^n - b\|_\infty).$$

(ii) *When $\rho \neq 1$, assume that $\varphi$ belongs to $\mathrm{H}^1(G)$. Then, for any finite open interval $(\ell', r') \subset G$, there exists a constant $C$ depending on $\lambda$, $\Lambda$, $\|\varphi\|_{\mathrm{H}^1}$, $\ell'$, $r'$ and $T$ such that*

$$\sup_{(t,x) \in [0,T] \times (\ell',r')} |u^n(t,x) - u(t,x)| \leq C(\|\rho^n - \rho\|_\infty + \|a^n - a\|_\infty + \|b^n - b\|_\infty).$$

PROOF. We set $G' = (\ell', r')$. For $T > 0$, we denote by $|\cdot|_{G',T}$ the norm defined by

$$|v|_{G',T} = \left(\sup_{t \in [0,T]} \|v(t,\cdot)\|_{\mathrm{L}^2(G')}^2 + \int_0^T \|\nabla v(t,\cdot)\|_{\mathrm{L}^2(G')}^2 \, dt\right)^{1/2}$$

for any $v \in \mathcal{C}(0,T; \mathrm{L}^2(G)') \cap \mathrm{L}^2(0,T; \mathrm{H}^1(G)')$. The convergence relies on the following estimate (which is easily derived from [19], Chapter II.3 inequality (3.1), page 74 and the Poincaré inequality, e.g.):

$$\sup_{(t,x) \in \mathbb{R}_+ \times G'} |v(t,x)| \leq C|v|_{G',T} \tag{31}$$

for a constant $C$ that depends only on $T$ and $G'$.

We now apply (31) to $v^n = u^n - u$. Indeed, it is easily seen that $v^n$ is the weak solution to

$$\begin{aligned}
\frac{\partial v^n(t,x)}{\partial t} &= L^n v^n(t,x) + \frac{1}{2}\nabla((a^n(x) - a(x))\nabla u(t,x)) \\
&\quad + (b^n(x) - b(x))\nabla u(t,x) + \left(\frac{1}{\rho^n} - \frac{1}{\rho}\right)\partial_t u(t,x),
\end{aligned} \tag{32}$$

where $\partial_t u(t,x)$ is a priori a distribution in $\mathrm{L}^2(0,T; \mathrm{H}^{-1}(G))$, where $\mathrm{L}^2(0,T; \mathrm{X})$ is the space of $\mathrm{L}^2$-functions with values in a Banach space X.

(i) In this case, the last term in the right-hand side of (32) vanishes. If $G$ is finite, we use $v^n$ as a test function in (32). After integrating with respect to $t$, we see that we are in the same position as in the elliptic case and standard computations yield the desired result.

If $G$ is not finite, we fix $\ell' < r'$ and we choose a smooth function $\xi$ with compact support such that $\xi(x) = 1$ on $(\ell', r')$. Then, since any function $v$ in $\mathrm{H}^1(G)$ is locally bounded, $\nabla(v\xi) = \xi\nabla v + v\nabla\xi$. We use $v^n \xi$ as a test function in (32) and the expansion of $\nabla(v\xi)$: this reduces the problem to the case where $G$ is finite. Hence, the proposition is proved.

(ii) As in [19], Theorem III.6.1, page 178, we show that $\partial_t u(t,x)$ belongs indeed to $\mathrm{L}^2(0,T; \mathrm{L}^2(G))$. For that, we assume at first that $a$, $\rho$ and $b$ are smooth, so that $u$ is also smooth. We transform $L$ into a divergence form operator with characteristic $(\widehat{a}, \widehat{\rho}, 0) = (e^{\widehat{\Psi}} a, e^{-\widehat{\Psi}} \rho, 0)$ as in Section 2. We assume that $\varphi \in \mathcal{C}^1_c(G; \mathbb{R})$. In this case, using $\partial_t u$ as a test function for the PDE $\partial_t u = Lu$, and integrating by parts against $e^{\widehat{\Psi}}/\rho$ with respect to $x$ and with respect to $t$, one obtains

$$\begin{aligned}
\int_0^T \int_G &\widehat{\rho}(x)|\partial_t u(t,x)|^2 \, dx \, dt \\
&= -\int_0^T \int_G \widehat{a}(x)|\partial_{t,x} u(t,x)|^2 \, dx \, dt \\
&\quad + \int_G \widehat{a}(x)|\varphi'(x)|^2 \, dx - \int_G \widehat{a}(x)|\partial_x u(T,x)|^2 \, dx.
\end{aligned} \tag{33}$$

As $\widehat{a}$ and $\widehat{\rho}$ are bounded and positive, (33) yields

$$\int_0^T \|\partial_t u(t,\cdot)\|^2_{\mathrm{L}^2(G)} \, dt \leq C\|\varphi\|^2_{\mathrm{H}^1(G)}, \tag{34}$$

where $C$ depends only on the constants $\lambda$ and $\Lambda$.

If the characteristic of $L$ is not smooth and $\varphi$ belongs only to $\mathrm{H}^1(G)$, an approximation procedure shows that (34) is still true.

Hence, with the additional assumption that the initial condition $\varphi$ belongs to $\mathrm{H}^1(G)$, it is possible to proceed as in (i).  □

**9. Numerical simulations: the algorithm.**  In this section we give two algorithms to simulate either $X_t$ for a given time $t$, or $(X_\tau, \tau)$ (or $X_t$ given by $t < \tau$), where $\tau$ is the first time the process $X$ reaches the points $r$ or $\ell$ where a Dirichlet b.c. holds.

9.1. *What is computed?.*  Our algorithm can be used to approximate the following quantities:

• For any initial distribution $\mu$ and any $t > 0$, we may compute the probability density function $u^*(t, y) = \int_G \mu(dx) p(t, x, y)$ with respect to $dy/\rho(y)$, since $u^*(t, y)$ is the solution to the PDE

$$\frac{\partial u^*(t,y)}{\partial t} = L^* u^*(t,y) \quad \text{and} \quad u^*(t,y)\,dy \xrightarrow[t \to 0]{} \mu.$$

In this case, we use the approximation, for a given $\varepsilon > 0$,

$$u^*(t,y) \simeq \frac{\rho(y)}{n 2\varepsilon} \operatorname{Card}\{i = 1, \ldots, n | Z^i \in [y - \varepsilon, y + \varepsilon]\},$$

where $Z^i$ are $n$ independent realizations of $X_t$ for which $t > \tau$ with the initial measure $\mu$.

• For any fixed $x \in (\ell, r)$ and any fixed $t > 0$, we may compute the solution $u(t, x)$ of the parabolic PDE (11) for any initial condition $\varphi$. In this case, we use the approximation

$$u(t,x) \simeq \frac{1}{n} \sum_{k=1}^{n} \varphi(Z^i),$$

where the $Z^i$'s are defined as previously with the initial distribution $\mu = \delta_x$.

Indeed, it is possible to deal with nonhomogeneous Dirichlet, Neumann or Robin b.c. by possibly coupling our algorithm locally with another algorithm: see Section 9.5.

• For any fixed $x \in (\ell, r)$, we may compute the solution $u(x)$ to the elliptic PDE

$$Lu = 0 \quad \text{on } (\ell, r), \quad u(r) = u_r \quad \text{and} \quad u(\ell) = u_\ell$$

for any $(u_r, u_\ell) \in \mathbb{R}^2$. In this case, we use the approximation

$$u(x) \simeq \frac{1}{n} \sum_{k=1}^{n} (u_r(Y^i)\mathbf{1}_{\{Y^i = u_r\}} + u_\ell(Y^i)\mathbf{1}_{\{Y^i = u_\ell\}}),$$

where the $Y^i$'s are $n$ independent realizations of $X_\tau$ under $\mathbb{P}_x$. If there is a Dirichlet b.c. at one of the endpoints, one may also consider a Neumann or Robin condition at the other endpoint, since a probabilistic representation of the solution similar to the one given in the parabolic case still holds.

Besides, if $\rho = 1$ and $b = 0$, then $L$ is self-adjoint with respect to the Lebesgue measure, and computing $\int_r^\ell \varphi(x) p(t, x, y) \, dx$ is sufficient to compute the solution $u(t, x)$ for any $x \in (\ell, r)$, since $p(t, x, y) = p(t, y, x)$ for all $t > 0$ and all $x, y \in (r, \ell)$.

REMARK 5. It is possible to get an analytical expression of $p(t, x, y)$ when $(x, y)$ belongs to a certain subspace of $G$ (see [15]), but it leads to rather complicated expressions involving spectral decompositions.

If $(a, \rho)$ are piecewise constant and $b = 0$, then the numerical errors come from the following: (1) The approximations of series by keeping a finite number of terms (see Section 9.3); (2) The fact that we use a pseudo-random generator instead of independent random variables; (3) The Monte Carlo error, that is, the replacement of $\mathbb{E}[Y]$ by $N^{-1} \sum_{k=1}^N Y^i$, where $Y$ is a random variable, and the $Y^i$'s are copies of independent random variables with the same distribution as $Y$.

The errors of type (1) and (2) are very small, and the error of type (3) occurs in all Monte Carlo methods and depends on the variance of $Y$.

9.2. *The algorithm.* For a general $(a, \rho, b)$, we may transform the diffusion process into a diffusion process with characteristics $(\widetilde{a}, \widetilde{\rho}, 0)$ and then use piecewise constant approximations $(\widetilde{a}^n, \widetilde{\rho}^n, 0)$ of $(\widetilde{a}, \widetilde{\rho}, 0)$. An additional error comes from these approximations, which is discussed in Section 8.

Let $X^n$ be the process generated by $(\widetilde{a}^n, \widetilde{\rho}^n, 0)$. We have seen in Section 7 that one can find a deterministic bijection $\Phi^n$ such that $Y_t^n = \Phi^n(X_t^n)$ is the solution to the SDE

$$Y_t^n = Y_0^n + B_t + \sum_{k \geq 0} \beta_k^n L_t^{y_k^n}(Y^n),$$

with $y_k^n = \Phi^n(x_k^n)$ and $\beta_k^n$ given by (26).

In this section we give an algorithm to simulate $(\tau \wedge T, Y_{\tau \wedge T}^n)$ for a given time $T > 0$, where $\tau$ is the first time the process $Y^n$ reaches a point at which there is a Dirichlet b.c. The computation of $(\tau \wedge T, X_{\tau \wedge T}^n)$ is done by setting $X_{\tau \wedge T}^n = \Phi^{-1}(Y_{\tau \wedge T}^n)$.

This algorithm is easily simplified to simulate $(Y_\tau^n, \tau)$.

Around each point $y_k^n$, $Y^n$ behaves like a skew Brownian motion of parameter $\beta_k^n$. The justification of our algorithm relies on Lemmas 1 and 2.

For each $y_k^n \in \Phi^n(\mathcal{J}^n)$, one picks two points $y_k^{n,-}$ and $y_k^{n,+}$ such that

$$y_{k-1}^n \leq y_k^{n,-} < y_k^n < y_k^{n,+} \leq y_{k+1}^n \quad \text{and} \quad y_k^{n,+} - y_k^n = y_k^n - y_k^{n,-}.$$

Let us denote by $\mathcal{K}^n$ the set $\mathcal{K}^n = \bigcup_{k \geq 0} \{y_k^{n,-}\} \cup \bigcup_{k \geq 0} \{y_k^{n,+}\}$. The sets $\mathcal{K}^n$ and $\Phi^n(\mathcal{J}^n)$ may have commons elements. The points of $\mathcal{K}^n$ have been introduced in order to use Lemma 1, and then to use random variables which are easily simulated.

The idea is to compute the successive times and positions on $\mathcal{K}^n \cup \Phi^n(\mathcal{J}^n)$ of a particle, with a special treatment to deal with the final time. The successive positions of the the particle on $\mathcal{K}^n \cup \Phi^n(\mathcal{J}^n)$ give simply a Markov chain on a discrete space, but dealing with the time requires a special treatment.

Let us fix $n$ large enough to get a fine approximation of the diffusion process. Since $n$ is fixed once and for all, we omit future reference to this integer.

The following algorithm takes as an input a horizon time $T$ and the starting point $\widetilde{y}$ of the particle. It returns $(\tau \wedge T, Y_{T \wedge \tau})$, where $\tau$ is the first time the particle hits one of the endpoints of $[\ell, r]$ where the Dirichlet b.c. holds.

NOTATION. We denote by $B$ a Brownian motion and $Z^\beta$ a skew Brownian motion of parameter $\beta \in [-1, 1]$. All the random variables simulated in this algorithm are assumed to be independent.

THE MAIN LOOP. We now explain how to update the position of the particle when it lies at a point $y = y_k$ at time $t$. We start with $y = \widetilde{y}$ and $t = 0$, and this loop is executed until the algorithm stops.

*Case $y \in \Phi(\mathcal{J})$*: Here, $y_k$ corresponds to a point in which the coefficients $a$ and $\rho$ may be discontinuous.

- If $y_k$ is at one of the endpoints of $[\ell, r]$ where the Dirichlet b.c. holds, then return $(t, y_k)$ and stop.
- Let $Z^{(\beta_k+1)/2}$ be a skew Brownian motion of parameter $(\beta_k+1)/2$ such that $Z_0^{(\beta_k+1)/2} = y_k$. Set $\tau = \inf\{s \geq 0 | Z_s^{(\beta_k+1)/2} \in \{y_k^-, y_k^+\}\}$ and compute $\gamma = \mathbb{P}_{y_k}[T - t < \tau]$. (Note that from Lemma 1, $\gamma$ does not depend on $\beta_k$ since $y_k - y_k^- = y_k^+ - y_k$.)
- Use a Bernoulli random variable of parameter $\gamma$ to decide if $\tau > T - t$ or $\tau < T - t$ (note that $\tau$ has not yet been simulated).
- If $\tau > T - t$, then the particle does not exit from $[y_k^-, y_k^+]$ before $T$: draw a realization $z$ of $Z_{T-t}^{(\beta_k+1)/2}$ given $\{\tau > T - t\}$. Finally, return $(T, z)$.
- If $\tau < T - t$, draw a realization $(t', z)$ of $(\tau, Z_\tau^{(\beta_k+1)/2})$ given $\{\tau < T - t\}$. Indeed, $\tau$ and $Z_\tau^{(\beta_k+1)/2}$ are independent, and $Z_\tau^{(\beta_k+1)/2}$ is a Bernoulli random variable of parameter $(\beta_k+1)/2$ with values in $\{y_k^-, y_k^+\}$. Update the current position of the particle by setting $y \leftarrow z$ and the current time by setting $t \leftarrow t + t'$. Restart at the beginning of the loop.

*Case $y \notin \Phi(\mathcal{J})$*: Except maybe for the initial point $\widetilde{y}$, this means that $y \in \mathcal{K} \setminus \Phi(\mathcal{J})$. Indeed, the operations are similar to the previous ones, except that one uses a Brownian motion instead of a skew Brownian motion.

- Compute $\gamma = \mathbb{P}_y[T - t < \tau]$ with $\tau = \inf\{s \geq 0 | B_s \in \Phi(\mathcal{J})\}$.
- Use a Bernoulli random variable of parameter $\gamma$ to decide if $\tau < T - t$ or $\tau > T - t$.
- If $\tau > T - t$, then draw a realization $z$ of $B_{T-t}$ given $T - t < \tau$ under $\mathbb{P}_y$. Then return $(T, z)$ and stop.
- If $\tau < T - t$, then draw a realization $z$ of $B_\tau$ given $\tau < T - t$ and afterward a realization $t'$ of $\tau$ given $\tau < T - t$ and $B_\tau = z$. Update the current position of the particle by setting $y \leftarrow z$ and the current time by setting $t \leftarrow t + t'$. Restart at the beginning of the loop.

9.3. *The random variables that should be simulated.* Let $(B, (\mathbb{P}_x)_{x \in \mathbb{R}})$ be a Brownian motion. Let $\tau_{a,b} = \inf\{t \geq 0 | B_\tau \in \{a, b\}\}$ for $a < b$. Using the scaling property, $\mathbb{P}_x \circ \tau_{a,b}^{-1}$ is equal to $\mathbb{P}_{(x-b)/(b-a)} \circ ((b-a)^2 \tau)^{-1}$. Hence, we set $\tau = \tau_{-1,1}$ and we assume that $a = -1$, $b = 1$.

Mainly, we need to simulate $\tau$ or some random variables whose distributions are related to that of $\tau$. From a numerical point of view, in order to simulate a random variable with distribution function $F$, we may compute $F^{-1}(U)$, where $U$ is a realization of a uniform random variable over $[0, 1]$ (about the simulation of random variables, see, e.g., [10]).

We give then explicit expressions of the densities of the random variables of interest (see [6] or [28], e.g.). Their distribution functions are easily computed, and may then be efficiently inverted using Newton's method.

If $G(t, x) = \mathbb{P}_x[\tau < t]$, then

$$\text{(35)} \quad \frac{\partial G}{\partial t}(t, x) = \sum_{k=-\infty}^{\infty} \left( \frac{1 + x + 4k}{\sqrt{2\pi} t^{3/2}} e^{-(1+x+4k)^2/2t} + \frac{1 - x + 4k}{\sqrt{2\pi} t^{3/2}} e^{-(1-x+4k)^2/2t} \right)$$

or

$$\text{(36)} \quad \frac{\partial G}{\partial t}(t, x) = \frac{\pi}{2} \sum_{k=0}^{+\infty} (-1)^k (2k+1) \exp\left(\frac{-\pi^2 (2k+1)^2 t}{8}\right) \cos\left(x\pi\left(k + \frac{1}{2}\right)\right).$$

Besides, if $H(t, x) = \mathbb{P}_x[\tau < t | B_\tau = 1]$, then

$$\text{(37)} \quad \frac{\partial H}{\partial t}(t, x) = \frac{2}{1+x} \sum_{k=-\infty}^{\infty} \left( \frac{1 - x + 4k}{\sqrt{2\pi} t^{3/2}} e^{-(1-x+4k)^2/2t} \right)$$

or

$$\text{(38)} \quad \frac{\partial H}{\partial t}(t, x) = \frac{-\pi}{2(1+x)} \sum_{k=1}^{+\infty} (-1)^k k \exp\left(\frac{-\pi^2 k^2 t}{8}\right) \sin\left(\frac{k\pi}{2}(x+1)\right).$$

The Bayes formula allows then to compute $\mathbb{P}_x[B_\tau = 1 | \tau < t]$ and by symmetry, $\mathbb{P}_x[B_\tau = -1 | \tau < t] = \mathbb{P}_{1-x}[B_\tau = 1 | \tau < t]$.

We also need $\mathbb{P}_x[B_t < y | t < \tau] = F(t, x, y) / \mathbb{P}_x[t < \tau]$, where $F(t, x, y)$ is defined as $F(t, x, y) = \mathbb{P}_x[B_t < y; t < \tau]$. The density of the Brownian motion killed when exiting from $[-1, 1]$ is

$$
\text{(39)} \quad \frac{\partial F}{\partial y}(t, x, y) = \frac{1}{\sqrt{2\pi t}} \sum_{k=-\infty}^{+\infty} \left( \exp\left(-\frac{(x - y - 4k)^2}{2t}\right) - \exp\left(-\frac{(x + y + 2 + 4k)^2}{2t}\right) \right)
$$

or, using a spectral decomposition,

$$
\text{(40)} \quad \frac{\partial F}{\partial y}(t, x, y) = \frac{1}{2} \sum_{k=1}^{+\infty} \exp\left(-\frac{k^2 \pi^2 t}{8}\right) \sin\left(\frac{k\pi}{2}(x+1)\right) \sin\left(\frac{k\pi}{2}(y+1)\right).
$$

The series giving $F$, $G$ and $H$ converge very quickly and do not create numerical problems. The series (35), (37) and (39) are numerically suitable for small time, while (36), (38) and (40) are numerically suitable for large time.

9.4. *Efficiency.* Most of the computation time is spent in the simulation of the exit time from an interval for the Brownian motion.

If we denote by $\mathcal{M}^n = \mathcal{J}^n \cup \mathcal{K}^n$, where $\mathcal{K}^n$ has been introduced at the end of Section 9.1, then we need to simulate the exit time of intervals of type $I_k^n = [z_{k-1}^n, z_{k+1}^n]$ if the $z_i^n$'s are the ordered points of $\mathcal{M}^n$. Thus, the "cost" of our algorithm may be identified with the number of times one needs to simulate the exit time of some $I_k^n$.

However, it is difficult to estimate this number $N^\tau$ of computations since it depends strongly on (a) the distance between $z_{k-1}^n$ and $z_{k+1}^n$ for $k \geq 1$, and (b) the number of passages on each interval $I_k^n = [z_{k-1}^n, z_{k+1}^n]$ for $k \geq 1$.

Besides, the number of such intervals $I_k^n$ depends on the variation of the coefficients. Intuitively, the flatter the coefficients, the more efficient is our algorithm.

In a simple but realistic case, we can give a rough estimate of the cost of our algorithm, which is matched by some numerical experiment (see Figure 3 in Section 10.2).

We assume that the coefficients $a$, $\rho$ and $b$ are of class $\mathcal{C}^1$ on $\mathbb{R} \setminus \{0\}$. Fix $n$ large enough and set $\Delta = 1/n$. We choose the map $\Phi^n$ so that the set $\Phi^n(\mathcal{J}^n) = \{k\Delta | k \in \mathbb{Z}\}$. We set $(a^n, \rho^n, b^n)(x) = (a, \rho, b)(k\Delta+)$ if $x \in [k\Delta, (k+1)\Delta]$ with $k \in \mathbb{Z}$. Let $u$ (resp. $u^n$) be the solution of the parabolic PDE $\partial u / \partial t = L u$ with $u(0, x) = f(x)$ [resp. $\partial u^n / \partial t = L^n u^n$ with $u^n(0, x) = \varphi(x)$], where $L$ (resp. $L^n$) is a differential operator with characteristics $(a, \rho, b)$ [resp. $(a^n, \rho^n, b^n)$]. Then, for any $(t,x)$, $|u^n(t,x) - u(t,x)|$ is of order $O(\Delta)$, since the coefficients are Lipschitz continuous on $\mathbb{R}_+^*$ and $\mathbb{R}_-^*$. Moreover, the Monte Carlo error in the evaluation of $u^n(t,x)$ is of order $O(1/\sqrt{N})$, where $N$ is the number of random independent particles, which gives a total error in the evaluation of $u(t,x)$ of order $O(1/\sqrt{N}) + O(\Delta)$.

As the intervals $I_k^n$ are of type $[(k-1)\Delta, (k+1)\Delta]$, in order to simulate $X_1^n$, we have to simulate $N^\tau$ independent random variables $\tau^1, \tau^2, \ldots$ giving the exit time from $[-\Delta, \Delta]$ for the Brownian motion starting from 0, where $N^\tau$ is such that $\tau^1 + \cdots + \tau^{N^\tau - 1} \leq 1 < \tau^1 + \cdots + \tau^{N^\tau}$. Replacing $\tau$ by its average $\mathbb{E}[\tau] = \Delta^2$, one gets that $N^\tau \sim 1/\Delta^2$. This means that the cost is of order $O(1/\Delta^2)$ for a trajectory. Thus, if one wants $|u(t,x) - u^n(t,x)|$ to be of order $O(\Delta)$, then one has to choose $N = 1/\Delta^2$. This means that, for a weak error of order $\Delta$, the number of times the random variable $\tau$ is simulated is of order $1/\Delta^4$ with our algorithm.

9.5. *Coupling with other algorithms.* Of course, this algorithm can be coupled with other algorithms (such as the Euler scheme, e.g.) if the coefficients are smooth enough outside some subset $I$ of $(\ell, r)$. Let us explain our idea: assume that the coefficients are smooth outside some interval $(\ell', r')$. Then, one can pick some points $r''$ and $\ell''$ such that $r < r'' < r'$ and $\ell' < \ell'' < \ell$. One can then use the algorithm which is the most adapted to the situation for the simulation $X$ outside $(\ell', r')$ and this shall be done until $X$ hits $r'$ or $\ell'$. Inside $(\ell', r')$ one can use the algorithm proposed here for the simulation $X$ until it hits $r''$ or $\ell''$.

One can also use a scheme where the process is killed at the discontinuities (there are efficient adaptations of the Euler scheme in this case). Thus, the distributions related to the skew Brownian motion may be used to re-inject the particle in the media.

Another way of coupling consists in using locally a scheme to deal with non-homogeneous Dirichlet, Neumann or Robin boundary conditions. For that, we use the following representations when the lateral function $\psi : \mathbb{R}_+^* \to \mathbb{R}$ is bounded and continuous (to simplify the notation, we assume that $r = +\infty$ and we set $\tau = \inf\{t > 0 | X_t = \ell\}$):

$$u(t,x) = \mathbb{E}_x[\varphi(X_t); t < \tau] + \mathbb{E}_x[\psi(\tau)] \quad \text{if } u(t,\ell) = \psi(t),$$

$$u(t,x) = \mathbb{E}_x[\varphi(X_t)] + \mathbb{E}_x\left[\int_0^t \psi(s)\, dL_s^\ell(X)\right] \quad \text{if } u'(t,\ell) = \psi(t),$$

$$u(t,x) = \mathbb{E}_x[e^{-\alpha L_t^\ell(X)} \varphi(X_t)] + \mathbb{E}_x\left[\int_0^t e^{-\alpha L_s^\ell(X)} \psi(s)\, dL_s^\ell(X)\right]$$

$$\text{if } \alpha u(t,r) + u'(t,r) = \psi(t) \text{ and } \alpha > 0.$$

With a Neumann or Robin b.c. at $\ell$, one has to consider the diffusion process which is reflected at this point. The Lépingle scheme [24, 25] gives an easy way to simulate the couple $(|B_t - \ell|, L_t^\ell(B))$ for any $t > 0$, where $B$ is a Brownian motion. Thus, one may then simulate $(X_{t+\delta t}, L_{t+\delta t}^\ell(X))$ by this way when $X_t = \ell$ and $\delta t$ is small enough. This gives approximations of integrals of type $\int_0^t \psi(s) \, dL_s^\ell(X)$, and allows to compute $u(t,x)$ using the previous formulae.

9.6. *Localization.* If $G = \mathbb{R}$, it follows from Aronson's estimates (see Proposition 2) that

$$\sup_{x \in \mathbb{R}} \mathbb{P}_x \left[ \sup_{s \in [0,t]} |X_s - x| \geq R \right] \leq C_1 \exp(-C_2 R^2/t)$$

for some constants $C_1, C_2$ depending only on $\lambda, \Lambda$ (see, e.g., [37] for a proof). Thus, one can assume that the coefficients $a$, $\rho$ and $b$ are constant far enough from the starting point, or that the process can be killed when it reaches the edges of a finite interval rather than dealing with a process that lives on the whole space $\mathbb{R}$.

## 10. Examples.

10.1. *The doubly skew Brownian motion.* We apply our algorithm to simulate the density of the solution to the SDE

$$X_t = X_0 + B_t + \tfrac{2}{3} L_t^{-1/2}(X) - \tfrac{2}{3} L_t^{1/2}(X),$$

where $B$ is a Brownian motion. Such a process may be called a *doubly skew Brownian motion*.

This process is then generated by

$$L = \frac{\rho}{2} \frac{d}{dx}\left(a \frac{d}{dx}\right)$$

with

$$(a(x), \rho(x)) = \begin{cases} (1,1), & \text{if } x \notin (-1/2, 1/2), \\ (2, 1/2), & \text{if } x \in (-1/2, 1/2). \end{cases}$$

In Figure 1 we represent the density at four different times. As expected, the density is more concentrated on the interval $(-1/2, 1/2)$, which agrees with our intuition, thanks to the choice of the coefficients.

10.2. *Diffusion with a coefficient discontinuous at one point.* We consider now that $b = 0$, $\rho = 1$ and that

$$a(x) = \begin{cases} 2 + \sin(x), & \text{if } x < 0, \\ 5 + \sin(x + \pi), & \text{if } x \geq 0. \end{cases}$$

We choose the points of $\mathcal{I}^n$ such that $\Phi^n(\mathcal{I}^n) = \{k\Delta \,|\, k \in \mathbb{Z}\}$ with $\Delta = 1/n$. The density $p(t, x, y)$ for $t = 1.0$ and $x = 0.5$ of the process is represented

in Figure 2 for $\Delta = 0.1$ and 10.000 particles, while the number of times the exit time from an interval for the Brownian motion have been drawn is shown in Figure 3. Figure 3 is fully in agreement with the rough estimate of Section 9.4.

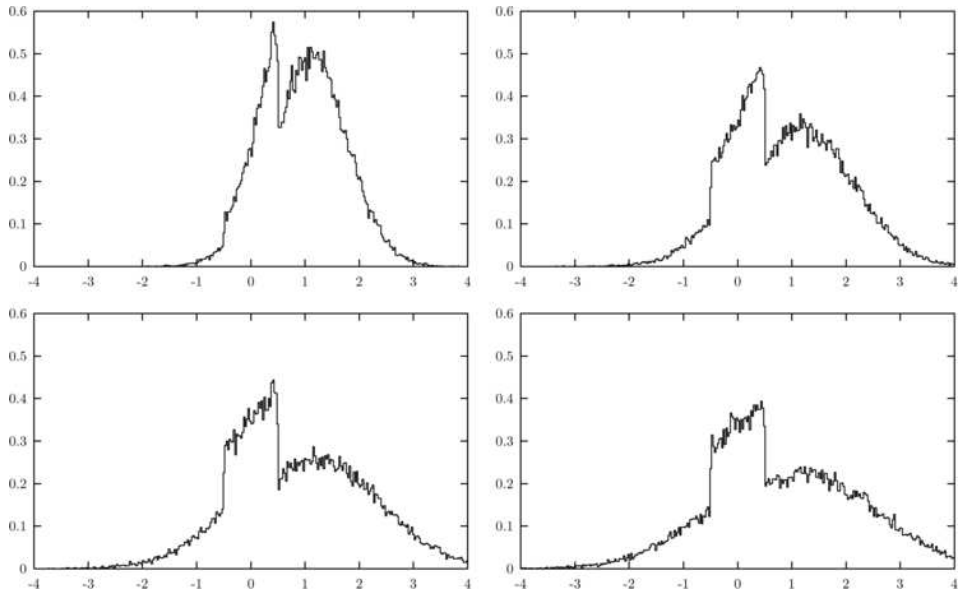

FIG. 1. *Density $p(t,x,y)$ of the doubly skew Brownian motion with $x = 1$ and $t = 0.5$, $t = 1$, $t = 1.5$ and $t = 2$ (with 50.000 particles).*

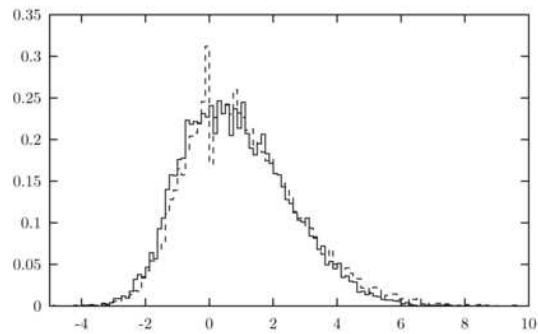

FIG. 2. *Densities $p(t,x,y)$ of the diffusion of Section 10.2 with $t = 1.0$ and $x = 0.5$: The full line is for the histogram obtained with the scheme presented in this article. The dashed line is for the histogram obtained with an Euler scheme.*

10.2.1. *Comparison with the Euler scheme.* In Figure 2 we have also drawn the density obtained with the Euler scheme presented by one of the authors in [27]. Indeed, the process $X$ is solution to the SDE

$$X_t = X_0 + \int_0^t \sigma(X_s)\,dB_s + \frac{1}{2}\int_0^t a'(X_s)\,ds + \frac{a(0+) - a(0-)}{a(0+) + a(0-)} L_t^0(X).$$

Using the function

$$\phi(x) = \frac{1-2\beta}{1-\beta}\mathbf{1}_{\mathbb{R}_+}(x) + \frac{1}{1-\beta}\mathbf{1}_{\mathbb{R}_-}(x)$$

with

$$\beta = (a(0+) - a(0-))/2a(0+),$$

one obtains from Itô–Tanaka that $Y_t = \phi(X_t)$ is the solution to some SDE with coefficients that are discontinuous at 0, but *without* local time, for which an Euler scheme is possible and may be applied. The convergence of this scheme with discontinuous coefficients is proved in [39], and [27] provides an estimation of the speed of convergence of this scheme in this particular case. In Figure 2 we use the time step $\delta t = 0.01$ and we see that the two empirical densities agree.

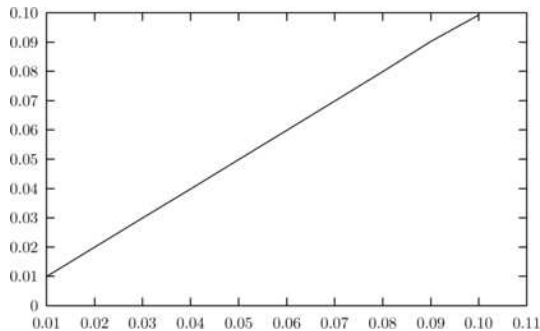

Fig. 3. $1/\sqrt{N}$ *in function of* $\Delta$, *where* $N$ *is the average number of simulations of an exit time.*

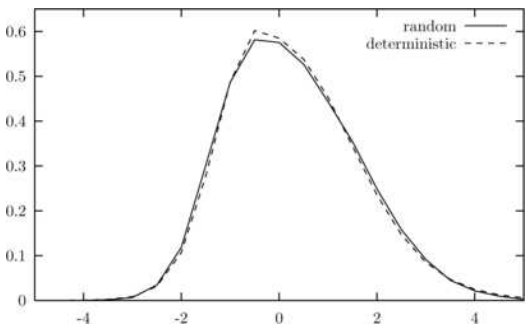

Fig. 4. *Interpolation of the solution* $u(t,x)$ *of* (41) *computed at points* $k/2$ *with* $k \in \mathbb{Z}$ *for* $t = 0.5$.

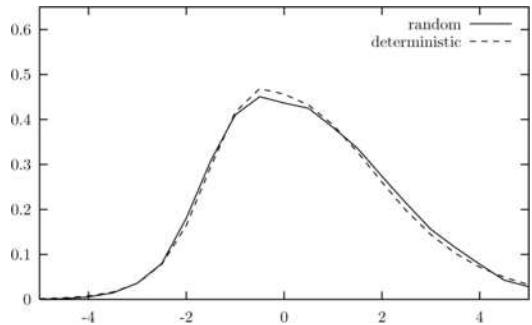

FIG. 5. *Interpolation of the solution $u(t,x)$ of* (41) *computed at points $k/2$ with $k \in \mathbb{Z}$ for $t = 1.0$.*

10.2.2. *Comparison with a deterministic scheme.* We consider now the PDE

$$\text{(41)} \qquad \frac{\partial u}{\partial t}(t,x) = Lu(t,x) \qquad \text{with } u(0,x) = \varphi(x),$$

where $\varphi(x) = \cos(x)$ if $|x| \leq \pi/2$ and $\varphi(x) = 0$ otherwise.

With our scheme, we computed $u(t,x)$ at time $t = 0.5$ and time $t = 1.0$ and at the points $k/2$ with $k \in \mathbb{Z}$. We also use the one-dimensional solver pdepe provided with MATLAB to solve (41) (indeed, we use Dirichlet boundary conditions at $-15$ and $15$, but this does not affect the computations for small times). In Figures 4 and 5, we see that the interpolated curves agree, even with a small amount of particles (here, 5.000 were used for each starting point).

Projet OMEGA
Institut Élie Cartan
Campus scientifique
BP 239
54506 Vandœuvre-lès-Nancy cedex
France
e-mail: Antoine.Lejay@iecn.u-nancy.fr

Projet OMEGA
INRIA Sophia-Antipolis
2004 Route des Lucioles
BP 93
06902 Sophia Antipolis
France
e-mail: Miguel.Martinez@sophia.inria.fr